\documentclass[11pt]{article}
\usepackage{amssymb}

\newtheorem{thm}{Theorem}

\newtheorem{lemma}[thm]{Lemma}

\newtheorem{cor}[thm]{Corollary}

\newtheorem{conj}[thm]{Conjecture}

\newcommand{\beq}[1]{\begin{equation}\label{#1}}
\newcommand{\enq}[0]{\end{equation}}

\newcommand{\qed}[0]{\begin{flushright} \rule{2mm}{3mm} \end{flushright}}

\newcommand{\C}[2]{{{#1}\choose{{#2}}}}
\newcommand{\ga}[0]{\alpha }
\newcommand{\gb}[0]{\beta }
\newcommand{\gc}[0]{\gamma }

\newcommand{\gl}[0]{\lambda }

\newcommand{\0}[0]{\emptyset}

\newcommand{\Rr}[0]{\mbox{${\bf R}$}}

\newcommand{\A}[0]{{\cal A}}
\newcommand{\B}[0]{{\cal B}}

\newcommand{\I}[0]{{\cal I}}
\newcommand{\m}[0]{{\cal M}}

\newcommand{\bn}[0]{\bigskip\noindent}
\newcommand{\mn}[0]{\medskip\noindent}

\newcommand{\sub}[0]{\subseteq}
\newcommand{\sm}[0]{\setminus}
\renewcommand{\dots}[0]{,\ldots,}

\newcommand{\uone}[0]{\underline{1}}

\begin{document}

\renewcommand{\thefootnote}{\fnsymbol{footnote}}

\title{A strong log-concavity property for measures on Boolean algebras*}

\author{
J. Kahn and M. Neiman\\
{Rutgers University } \\
{\footnotesize email: jkahn@math.rutgers.edu; neiman@math.rutgers.edu}
}
\date{}

\footnotetext{ * Supported by NSF grant DMS0701175.}

\maketitle

\begin{abstract}
We introduce the {\em antipodal pairs property}
for probability measures on finite Boolean algebras
and prove that
conditional versions imply strong forms of log-concavity.
We give several applications of this fact,
including improvements of some results of Wagner \cite{Wagner2};
a new proof of a theorem of Liggett \cite{LiggettULC} stating that
{\em ultra-log-concavity} of sequences
is preserved by convolutions;
and some progress on a well-known log-concavity
conjecture of J. Mason \cite{Mason}.

\end{abstract}

\section{Introduction}\label{section:intro}

\mn
{\em Log-concavity and the antipodal pairs property}

Before stating our main result, we fix some notation and recall some definitions.
Given a finite set $S$,
denote by $\m=\m_S$ the set of probability measures
on $\Omega=\Omega_S=\{0,1\}^S$.
As a default we take
$S=[n] =\{1\dots n\}$ (which for us is simply
a generic $n$-set), using
$\Omega$ and $\m$ for $\Omega_{[n]}$ and $\m_{[n]}$.
We will occasionally identify
$\Omega$ with the Boolean algebra $2^{[n]}$
(the collection of subsets of $[n]$ ordered by inclusion)
in the natural way
(namely, identifying a set with its indicator).

We will be interested in several properties of measures which are preserved by the operation of {\em conditioning}, which for us always means fixing the values of some variables
(this specification always assumed to have
positive probability); thus a measure obtained from
$\mu\in \m$ by conditioning is one of the form
$\mu(\cdot|\eta_i=\xi_i ~\forall i\in I)$
(which we regard as a measure on $\Omega_{[n] \sm I}$)
for some $I\sub [n]$ and $\xi \in \{0,1\}^I$.
(If we think of $\Omega$ as $2^{[n]}$, then
conditioning amounts to restricting our measure to
some interval $[J,K]$ of $2^{[n]}$ (and normalizing).)

Recall that a sequence $a=(a_0\dots a_n)$ of real numbers
(here always nonnegative)
is {\em unimodal} if there is some $k\in \{0\dots n\}$ for which
$a_0 \leq a_1 \leq \cdots \leq a_k \geq \cdots \geq a_n $,
and {\em log-concave} (LC) if
$a_i^2\geq a_{i-1}a_{i+1}$ for $1\leq i\leq n-1$.
Of course a nonnegative LC sequence with no internal zeros is unimodal
(where ``no internal zeros" means $\{i : a_i \neq 0\}$ is an interval).
Following \cite{Pem} we say $a$ (as above)
is {\em ultra-log-concave} (ULC) if
the sequence $(a_i/\C{n}{i})_{i=0}^n$ is log-concave
and has no internal zeros.
We also say
$\mu\in\m$ is ULC if
its {\em rank sequence}, $(\mu(|\eta|=i))_{i=0}^n$, is ULC
(where $|\eta|=\sum\eta_i$).
We define ``$\mu$ is LC" and ``$\mu$ is unimodal" similarly,
except that for the former we add the
requirement
that
the rank sequence
have no internal zeros.

For $\mu\in \m$ set
\beq{APseq}
\ga_i(\mu) =
\C{n}{i}^{-1}
\sum\{\mu(\eta)\mu(\uone-\eta):\eta\in \Omega, |\eta|=i\}
\enq
(where $\uone=(1 \dots 1)$).
Say that $\mu\in \m_{2k}$ has the {\em antipodal pairs property}
(APP) if
$\ga_k(\mu)\geq \ga_{k-1}(\mu),$
and that $\mu\in \m$ has the
{\em conditional antipodal pairs property} (CAPP)
if any measure obtained from $\mu$ by conditioning on
the values of some $n-2k$ variables (for some $k$) has the APP.
Our main result is
\begin{thm} \label{CAPPULC}
For measures without internal zeros in their rank sequence,
the CAPP implies ULC.
\end{thm}
A somewhat more general version of
Theorem \ref{CAPPULC}
is stated and proved in
Section \ref{section:proof}.
In the
rest
of this introduction
we provide a little context
and sketch some consequences to be established in later sections.

\mn
{\em Negative dependence properties}

We need to
briefly review a few
negative dependence notions;
for much more on this see
e.g.
\cite{Pem,choewagner03,Wagner,BBL,KN1}.
Recall that events $\A,\B$ in a probability space
are \emph{negatively
correlated}---we write $\A\downarrow \B$---if $\Pr(\A\B)\leq\Pr(\A)\Pr(\B)$.
We say $\mu$ has {\em negative correlations} (or {\em is NC}) if
$\eta_i\downarrow \eta_j$
(that is, $\{\eta_i=1\}\downarrow \{\eta_j=1\}$)
whenever $i\neq j$.
A stronger property is obtained by requiring NC for every measure $W \circ \mu$ of the form
$$W\circ \mu(\eta) ~\propto~ \mu(\eta)\prod W_i^{\eta_i}$$
with $W=(W_1\dots W_n)\in\Rr_+^n$;
we say $W \circ \mu$ is obtained from $\mu$ by {\em imposing
an external field} and
(following \cite{choewagner03}, \cite{Wagner}) say $\mu$ has the
{\em Rayleigh property} if every measure gotten from $\mu$ by imposing
an external field is NC.

Theorem \ref{CAPPULC} was discovered during attempts
to prove a sequence of conjectures of
Pemantle stating that various
negative dependence properties, including the Rayleigh property,
imply ULC; see \cite[Conjecture 4]{Pem}.
That Rayleigh
implies ULC was also
conjectured by Wagner \cite{Wagner}.
As it happens, even the weakest version of Pemantle's conjecture is false;
see \cite{BBL}, \cite{KN1} for counterexamples and further discussion.
Still, Theorem \ref{CAPPULC} does turn out to be helpful
in establishing ULC in some other settings, which we now indicate.

\medskip
A first, easy consequence is improvement of some of the results
of \cite{Wagner2},
for which we need to recall some terminology from that paper.
Given a positive integer $k$ and positive real number $\gl$, say that $\mu$ satisfies
$\gl$-$\textrm{Ray}[k]$
if every measure $\nu$
gotten from $\mu$ by imposing an external field and then projecting
onto a set $S$ of $2k$ variables satisfies
\beq{lray}
\sum\{\nu(\eta)\nu(\uone-\eta):\eta\in \Omega, |\eta|=k\} \geq \gl \sum\{\nu(\eta)\nu(\uone-\eta):\eta\in \Omega, |\eta|=k-1\}.
\enq
(For the (standard) definition of projection, see
Section \ref{section:corproof}.)
With the notation of (\ref{APseq}), the above condition is
$$\ga_k(\nu) \geq \frac{\gl k}{k+1} \ga_{k-1}(\nu);$$
thus $(1+1/k)$-$\textrm{Ray}[k]$ says that
each $\nu$
as above has the APP.
(As observed by Wagner \cite{Wagner2}---see his 
Proposition 4.6---$\gl =(1+1/k)$ is an 
``especially natural strength for these conditions'';
there as here, this is essentially because $1+1/k$ is the 
ratio of the numbers of summands on the two sides of 
(\ref{lray}).)
Note also that $2$-$\textrm{Ray}[1]$ is precisely the Rayleigh property.

Say that $\mu\in \m$ is BLC$[m]$ if every measure gotten from $\mu$ by imposing an external field and then projecting onto a set of size at most $m$ is ULC
(the acronym is for ``binomial log-concavity''),
and BLC if it is BLC$[m]$ for all $m$.
In \cite{choewagner03} and \cite{KN1}, BLC$[m]$ is called LC$[m]$.)
Wagner proved
\begin{thm}[\cite{Wagner2}, Theorem 4.3]\label{ThmWagner}
If a measure satisfies $2$-$\textrm{Ray}[1]$ and $(1+1/k)^2$-$\textrm{Ray}[k]$ for all $2 \leq k \leq m$, then it is BLC$[2m+1]$.
\end{thm}
(In \cite{Wagner2}
this is stated only for uniform measure on the bases of a matroid, but the proof is valid for general measures in $\m$.)
Theorem \ref{CAPPULC} implies the following strengthening 
(see Section \ref{section:corproof}).
\begin{cor}\label{CorWagner}
If a measure satisfies $(1+1/k)$-$\textrm{Ray}[k]$ for all $1 \leq k \leq m$, then it is BLC$[2m+1]$.
\end{cor}
Using Corollary \ref{CorWagner} in place of Theorem \ref{ThmWagner} improves
Corollary 4.5(b) and Theorem 5.2 of \cite{Wagner2} by substituting BLC
for the weaker property $\sqrt{\textrm{BLC}}$; see \cite{Wagner2} for
definitions and statements.

\mn
{\em Ultra-log-concave sequences}

It is easy to see that if
$\mu \in \m_S$ and $\nu \in \m_T$ are Rayleigh then the product measure $\mu \times \nu$ (given by
$\mu \times \nu (\xi,\eta) = \mu(\xi) \nu(\eta)$ for
$(\xi,\eta) \in \{0,1\}^S \times \{0,1\}^T)$)
is also Rayleigh.
Note that the rank sequence of $\mu \times \nu$ is the convolution of the rank sequences for $\mu$ and $\nu$.
One consequence of the (false) conjecture mentioned earlier,
that Rayleigh measures are ULC, would have been that the convolution of two ULC sequences is ULC or, equivalently, that the product of two ULC measures is ULC.
(The implication
follows from a result of Pemantle \cite[Theorem 2.7]{Pem}
stating
that for {\em exchangeable} measures (those for which $\mu(\eta)$ depends
only on $|\eta|$) the properties Rayleigh and ULC coincide.)
Surprisingly---given that the analogous statement for ordinary log-concavity
is fairly trivial---preservation of ULC under convolution turns out
not to be so obvious; it was conjectured by
Pemantle \cite{Pem} (motivated by the preceding considerations)
and proved by Liggett:
\begin{thm}[\cite{LiggettULC}, Theorem 2] \label{ConvULC}
The convolution of two ULC sequences is ULC.
\end{thm}
In Section \ref{section:ConvULC} we derive this from
Theorem \ref{CAPPULC} and also
discuss a potentially interesting strengthening
of ULC for measures that is again preserved by products.

\mn
{\em Mason's conjecture}

For the purposes of this introduction 
we regard a matroid as a collection
$\I$ of independent subsets
of some ground set $E$.
We will not go into matroid definitions; see e.g.
\cite{Welshbook} or \cite{Oxley}.
Prototypes are the collection of (edge sets of)
forests of a graph (with edge set $E$)---this is a
{\em graphic} matroid---and (as it turns out, more generally)
the collection of linearly independent subsets of
some finite subset $E$ of some (not necessarily finite)
vector space; for present purposes not too much is lost
by thinking only of graphic matroids.

We are interested in the
{\em independence numbers} of
$\I$, that is, the numbers
$$a_k = a_k(\I) = |\{I\in\I:|I|=k\}| ~~~ k=0\dots n,$$
for which a celebrated conjecture of
J. Mason \cite{Mason} says

\begin{conj}\label{Cmason}
For any matroid $\I$ on a ground set of size n,
the sequence $a=a(\I)=(a_0\dots a_n)$ is ULC.
\end{conj}
(Note that $a$ will typically end with some 0's,
and also that in the graphic case $n$ counts edges, not vertices.)
Of course one can weaken Conjecture \ref{Cmason}
by asking for LC or unimodality
in place of ULC.
In fact Mason also stated the LC version, and unimodality,
first suggested by Welsh \cite{Welsh},
was the original conjecture
in this direction (and
even this, even for graphic matroids, remains open).

From the present viewpoint, Mason's Conjecture says that
uniform measure on $\I$ (regarded
in the usual way as a subset of $\{0,1\}^E$) is ULC.
(When $\I$ is graphic such a measure is a
{\em uniform spanning forest} (USF) measure
(``spanning" because we think of a member of $\I$ as a
subgraph that includes all vertices).)
In particular, according to Theorem \ref{CAPPULC}, Mason's conjecture
would follow from
\begin{conj}\label{matroidCAPPconj}
Uniform measure on the independent sets of a matroid has the CAPP.
\end{conj}
See also the remark following Corollary \ref{CAPUULC} for a possible
strengthening.

Of course here
it's enough to show APP, since each conditional measure
is just uniform measure on the independent sets of some minor.
Also, note that Mason's conjecture for graphic matroids would
have followed from the (false) conjecture of Pemantle \cite{Pem} and Wagner \cite{Wagner} that Rayleigh measures are ULC, if it could be shown that,
as conjectured in \cite{AN} (see \cite{Wagner,GW,KN1}
for more on this),
USF measures are Rayleigh.

Though probably not for lack of effort,
progress on Mason's conjecture has been fairly modest.
Dowling \cite{Dowling} proved that for each
$\I$ the sequence $(a_0 \dots a_8)$ is LC;
Mahoney \cite{Mahoney} proved that for graphic matroids
corresponding to {\em outerplanar} graphs,
the full sequence of independence numbers is LC;
and Hamidoune and Sala\"un \cite{HamidouneSalaun}
proved that for any matroid on a ground set of size $n$ the
sequence $(a_i/\C{n}{i})_{i=0}^4$ is LC, i.e. the sequence
$(a_i)$ is ``ULC up to 4''.

Here we adapt one of Dowling's arguments to prove Conjecture
\ref{matroidCAPPconj} for small matroids:
\begin{thm}\label{matroidCAPPsmall}
For every matroid on a ground set of size at most 11, uniform measure on  independent sets has the CAPP.
\end{thm}
This is proved in Section \ref{section:Mason}.
Combined with Theorem \ref{CAPPULC} (for (a))
or the more general Theorem \ref{CAPPULCLocal} below (for (b))
it gives
\begin{thm}\label{partialMason}
{\rm (a)} Every matroid on a ground set of size at most 11 satisfies
Conjecture \ref{Cmason}.
\\
{\rm (b)} For any matroid on a ground set of size $n$
with independence numbers $a_i$,
the sequence $(a_i/\C{n}{i})_{i=0}^6$ is LC (a.k.a. the sequence
$(a_i)$ is ``ULC up to 6'').
\end{thm}

\section{Proof of Theorem \ref{CAPPULC}} \label{section:proof}

We actually prove
a more general result that will be needed
in Section \ref{section:Mason}.
\begin{thm}\label{CAPPULCLocal}
Suppose $\mu \in \m$ has the property that, for every $k \in [t]$,
every measure gotten
from $\mu$
by conditioning on the values of $n-2k$ coordinates has the APP.
Then the sequences $(\mu(|\eta|=i)/\C{n}{i})_{i=0}^{t+1}$ and $(\mu(|\eta|=i)/\C{n}{i})_{i=n-t-1}^n$ are LC.
\end{thm}
For the rest
of this section it will be convenient to treat $\Omega$ as $2^{[n]}$,
so that (\ref{APseq}) becomes
$$\ga_i(\mu)=\C{n}{i}^{-1} \sum \{ \mu(X) \mu([n] \sm X) : X \in \C{[n]}{i} \}$$
(where $\C{[n]}{i}=\{X \sub [n]:|X|=i\}$).

We first need to recall some properties of the Johnson association scheme;
this material (up to (\ref{pij})) is taken from chapter 30 of \cite{vLW}.
Fix positive integers $n$ and $l$ with $l \leq n/2$, let
$\mathfrak{X}=\C{[n]}{l}$, and, for $i=0,1 \dots l$, let $A_i$ be the $\mathfrak{X} \times \mathfrak{X}$
adjacency matrix of $i$th associates, {\em viz.}
$$A_i(X,Y)=\left\{
\begin{array}{ll}
1 & \textrm{if } |X \cap Y|=l-i \\
0 & \textrm{otherwise.}
\end{array}
\right.$$
We write
elements of
$\Rr^{\mathfrak{X}}$ as row vectors.
For $T \sub [n]$ with $|T| \leq l$, let $e_T$ be the vector in $\Rr^{\mathfrak{X}}$ with
$$e_T(S)=\left\{
\begin{array}{ll}
1 & \textrm{if } S \supseteq T \\
0 & \textrm{otherwise,}
\end{array}
\right.$$
and let $U_i $ be the span of $\{e_T : T \in \C{[n]}{i}\} $.
Then $\dim{U_i}=\C{n}{i}$ and
$U_0 \subseteq U_1 \subseteq \cdots \subseteq U_l=\Rr^{\mathfrak{X}}$.
Set $V_0 = U_0$ and $V_i = U_i \cap U_{i-1}^{\bot}$ for $i=1,2 \dots l$, and
let $E_i$ be the projection of $\Rr^{\mathfrak{X}}$ onto $V_i$.
Then
$$\Rr^{\mathfrak{X}}=V_0 \oplus V_1 \oplus \cdots \oplus V_l$$
is an orthogonal decomposition,
$$E_i E_j=\left\{
\begin{array}{ll}
E_i & \textrm{if } i=j \\
0 & \textrm{if } i \neq j,
\end{array}
\right.$$
and
$$E_0 + E_1 + \cdots + E_l = I.$$
Note that $V_0$
consists of the constant vectors.
The span, $\mathfrak{A}$, of $A_0\dots A_l$
is an algebra under matrix multiplication (the {\em Bose-Mesner algebra}).
The set of matrices $\{E_0,E_1 \dots E_l\}$ is also a basis for $\mathfrak{A}$, with
\beq{AsEs}
A_i=\sum_{j=0}^l P_i(j)E_j  \ \ \ (i=0,1 \dots l),
\enq
where
\beq{pij}
P_i(j)=\sum_{k=0}^i(-1)^{i-k} \C{l-k}{i-k} \C{l-j}{k} \C{n-l+k-j}{k}.
\enq
The next lemma is presumably well-known.
\begin{lemma} \label{PSD}
For any $\gc_0, \gc_1 \dots \gc_l\in \Rr$, the $\mathfrak{X} \times \mathfrak{X}$ real symmetric matrix $M=\sum_{i=0}^l \gc_i A_i$ is positive
semidefinite if and only if $~\sum_{i=0}^l \gc_i P_i(j) \geq 0$ for
$j=0,1 \dots l$.
\end{lemma}

\mn {\em Proof.} Since $M=\sum_{j=0}^l \sum_{i=0}^l \gc_i P_i(j) E_j$ and $E_j$ is the orthogonal projection of $\Rr^{\mathfrak{X}}$ onto $V_j$, the eigenvalues of $M$ are $\{\sum_{i=0}^l \gc_i P_i(j) : j=0,1 \dots l\}$. \qed

\mn {\em Remark.} It is not hard to show,
using some additional properties of the Johnson scheme,
that the condition appearing in Lemma \ref{PSD} is equivalent to
the statement that
the vector $(\C{l}{i} \C{n-l}{l-i} \gc_i)_{i=0}^{l}$ satisfies
{\em Delsarte's inequalities} (\cite{Delsarte} or
\cite[p. 416]{vLW}).

\medskip
We also need one technical lemma:
\begin{lemma} \label{SUM}
For all positive integers $M,N$ and real numbers $a,b$,
$$\sum_{t=0}^N (-1)^t \frac{at+b}{t+M} \C{N}{t} = (\frac{b}{M}-a) \C{M+N}{M}^{-1}.$$
\end{lemma}

\mn {\em Proof.} It suffices to prove either of the equivalent
\beq{sum_ineq1}
\sum_{t=0}^N (-1)^t \frac{t}{t+M} \C{N}{t} \C{M+N}{M} = -1,
\enq
$$\sum_{t=0}^N (-1)^t \frac{M}{t+M} \C{N}{t} \C{M+N}{M} = 1,$$
since the desired identity is a linear combination of these.
We prove (\ref{sum_ineq1}), fixing
$M$ and proceeding by induction on $N$,
with the base case $N=1$ trivial.
For the induction step, we just check that the left
side of (\ref{sum_ineq1}) does not change when we replace $N$ by $N+1$;
indeed, the difference is
$$(-1)^{N+1} \frac{N+1}{M+N+1} \C{M+N+1}{M} + \sum_{t=1}^N (-1)^t \frac{t}{t+M}
\Big[ \C{N+1}{t} \C{M+N+1}{M} - \C{N}{t} \C{M+N}{M} \Big]$$
$$=(-1)^{N+1} \C{M+N}{M} +\sum_{t=1}^N (-1)^t \C{N}{t-1} \C{M+N}{M},$$
which is zero. \qed

\mn {\em Proof of Theorem \ref{CAPPULCLocal}.} Let $\mu$ be a measure on $2^{[n]}$ satisfying the hypotheses
of the theorem,
with rank sequence $(a_i)_{i=0}^n$.
Our goal is to show
\beq{ulc1}
l (n-l) a_l^2 \geq (l+1)(n-l+1) a_{l-1} a_{l+1}
\enq
for $l \in \{1 \dots t\} \cup \{n-t \dots n-1\}$; but, since
$\mu'\in \m$ given by
$\mu'(X)=\mu([n] \sm X)$ again satisfies the hypotheses
of Theorem \ref{CAPPULCLocal} and has rank sequence
$(a_{n-i})_{i=0}^n$,
it suffices to prove (\ref{ulc1}) when $l \leq \min\{t,n/2\}$.
To this end, fix such an $l$ and set
$$Z_{j,k}^i = \sum \{\mu(X) \mu(Y) : (X,Y) \in \C{[n]}{j} \times \C{[n]}{k} \textrm{ and } |X \cap Y| = i\};$$
with this notation, (\ref{ulc1}) is
\beq{ulc2}
l(n-l) \sum_{i=0}^l Z_{l,l}^i \geq (l+1)(n-l+1) \sum_{i=0}^{l-1} Z_{l-1,l+1}^i.
\enq

For each $i \in \{0,1 \dots l-1\}$ and $I \sub J \sub [n]$ with $|I|=i$ and $|J|=2l-i$, let $\mu_{I,J} \in \m_{J \sm I}$ be
the conditional measure with
$$\mu_{I,J}(X) \propto \mu(X \cup I) ~~~~ (X \sub J \sm I)$$
(or $\mu_{I,J}\equiv 0$ if $\mu([I,J])=0$).
By hypothesis (or trivially if $\mu_{I,J}\equiv 0$) 
$\mu_{I,J}$ has the APP, i.e.
\beq{condAPP1}
\ga_{l-i}(\mu_{I,J}) \geq \ga_{l-i-1}(\mu_{I,J}).
\enq
With
$$Z_{j,k}(I,J) = \sum \{\mu(X) \mu(Y) : (X,Y) \in \C{[n]}{j} \times \C{[n]}{k}, ~ X \cup Y=J, \textrm{ and } X \cap Y = I\},$$
we have
$$\ga_{j}(\mu_{I,J}) = \left( \sum_{I \sub X \sub J} \mu(X) \right)^{-1} \C{2l-2i}{j}^{-1} Z_{i+j,2l-i-j}(I,J),$$
and (\ref{condAPP1}) becomes
\beq{condAPP2}
\frac{l-i}{l-i+1} Z_{l,l}(I,J) \geq Z_{l-1,l+1}(I,J).
\enq
Summing
(\ref{condAPP2}) over $I, J$ with $|I|=i$ and $|J|=2l-i$ gives
\beq{ineq_from_capp}
\frac{l-i}{l-i+1} Z_{l,l}^i \geq Z_{l-1,l+1}^i
\enq
(since each pair $(X,Y)$ contributing to $Z_{l,l}^i$ contributes to the left
side of (\ref{condAPP2}) for exactly one choice of $(I,J)$, and similarly
for pairs contributing to $Z_{l-1,l+1}^i$).
Replacing each $Z_{l-1,l+1}^i$ in (\ref{ulc2})
by the (corresponding) left side of (\ref{ineq_from_capp}),
we find that it is enough to show
\beq{js_ineq}
\sum_{i=0}^l \gb_i Z_{l,l}^i \geq 0,
\enq
where
$$\gb_i=\frac{i(n+1)-l(l+1)}{l-i+1}.$$

In fact, we will show that
(\ref{js_ineq}) holds for {\em every} $\mu \in \m$.
Let $\psi=\psi_{\mu}$ be the vector in $\Rr^{\mathfrak{X}}$
with $\psi(X)=\mu(X)$ for $X \in \C{[n]}{l}$,
and recall the matrices $A_i$ defined before Lemma \ref{PSD}.
Since
$$\psi A_i \psi^T = Z_{l,l}^{l-i},$$
the left side of (\ref{js_ineq}) is
\beq{musummu}
\psi (\sum_{i=0}^l \gb_{l-i} A_i) \psi^T
\enq
and (\ref{js_ineq}) will follow from Lemma \ref{PSD} once we show
\beq{pij_ineqs}
\sum_{i=0}^l \gb_{l-i} P_i(j) \geq 0
\enq
for $j=0,1 \dots l$.

Fix $j\in [l]$.
(We deal with the case $j=0$ separately below.)
The left side of (\ref{pij_ineqs}) is
\beq{in1}
\sum_{i=0}^l \frac{(l-i)(n+1)-l(l+1)}{i+1} \sum_{k=0}^i(-1)^{i-k} \C{l-k}{i-k} \C{l-j}{k} \C{n-l+k-j}{k},
\enq
which we want to show is nonnegative.
Interchanging the order of summation and making the substitution $t=i-k$,
we may rewrite (\ref{in1}) as
\beq{in2}
\sum_{k=0}^l \C{l-j}{k} \C{n-l+k-j}{k} \sum_{t=0}^{l-k} (-1)^t \frac{(l-t-k)(n+1)-l(l+1)}{t+k+1} \C{l-k}{t}.
\enq
It is thus enough to show
\beq{in3}
\sum_{t=0}^{l-k} (-1)^t \frac{(l-t-k)(n+1)-l(l+1)}{t+k+1} \C{l-k}{t} \geq 0
\enq
whenever $k\leq l-1$
(since the $k=l$ term in (\ref{in2}) is zero).
But Lemma \ref{SUM}, with $N=l-k$, $M=k+1$, $a=-(n+1)$, and
$b=N(n+1)-l(l+1)$, says that the left side of (\ref{in3}) is
$$\frac{(n-l+1)(l+1)}{k+1} \C{l+1}{k+1}^{-1},$$
which is positive since $l \leq n/2$.

Finally we show that when $j=0$, (\ref{pij_ineqs}) holds with equality.
To see this, notice that when $\mu$ is uniform measure on $2^{[n]}$,
we have equality in
(\ref{ulc2}) and (\ref{ineq_from_capp}), and consequently (\ref{js_ineq}),
from which it follows (see (\ref{musummu})) that
\beq{mumu}
\psi (\sum_{i=0}^l \gb_{l-i} A_i) \psi^T=0.
\enq
But, since $\psi E_j\psi^T$ is $2^{-2n}\C{n}{l}$ if $j=0$ and zero otherwise,
the left side of (\ref{mumu}) is (by (\ref{AsEs}))
$$2^{-2n}\C{n}{l} \sum_{i=0}^l \gb_{l-i} P_i(0),$$ which gives the promised equality
in (\ref{pij_ineqs}).\qed

\section{Proof of Corollary \ref{CorWagner}} \label{section:corproof}

In this short section we use Theorem \ref{CAPPULC} to prove Corollary \ref{CorWagner}.
Recall that the {\em projection} of $\mu \in \m_n$ on $J\sub [n]$ is the measure
$\mu'$ on $\{0,1\}^J$ obtained by integrating out
the variables of $[n]\sm J$; that is,
$$\mu'(\xi) =
\sum\{\mu(\eta):\eta\in \Omega, \eta_i=\xi_i ~\forall i\in J\}
~~~~(\xi\in \{0,1\}^J).$$

\mn
{\em Proof of Corollary \ref{CorWagner}.}
The statement is:
if $\mu\in \m$ satisfies $(1+1/k)$-$\textrm{Ray}[k]$ for all $k\in [m]$,
$T\sub [n]$, $|T|\leq 2m+1$, and $\nu \in \m_T$ is obtained from $\mu$
by imposing an external field and projecting on $T$,
then $\nu $ is ULC.
By Theorem \ref{CAPPULC}, it suffices to show $\nu$ has the CAPP. 
(We should also show that $\nu$'s rank sequence has no internal zeros, but this follows immediately from \cite[Proposition 4.7(a)]{Wagner}.)
Any measure gotten from $\nu$ by conditioning on the values of the
variables in some set $T \sm S$
is the limit of a sequence
of measures, each
gotten from $\mu$ by imposing an external field and projecting
on $S$; the CAPP for $\nu$
thus follows from our assumption on $\mu$.\qed

\section{Convolution of ULC sequences} \label{section:ConvULC}

In this section we define a property of measures which
is stronger than ULC, prove it is preserved by products, and
show that this implies Theorem \ref{ConvULC}.

We begin with some definitions.
With $\ga_i(\mu)$ as in (\ref{APseq}), say $\mu \in \m$ is
{\em antipodal pairs unimodal} (APU) if the sequence $(\ga_i(\mu))_{i=0}^n$ is
unimodal (since $\ga_i(\mu)=\ga_{n-i}(\mu)$, this means
$\ga_0(\mu) \leq \cdots
\leq \ga_{\lfloor n/2 \rfloor}(\mu) = \ga_{\lceil n/2 \rceil}(\mu) \geq \cdots \geq \ga_n(\mu)$), and say $\mu$ is {\em conditionally antipodal pairs unimodal} (CAPU) if every measure obtained from $\mu$ by conditioning
is APU.
Since CAPU trivially implies the CAPP, Theorem \ref{CAPPULC} gives
\begin{cor}\label{CAPUULC}
Every CAPU measure is ULC.
\end{cor}
(As far as we know, Conjecture \ref{matroidCAPPconj} can be strengthened by replacing ``CAPP'' with ``CAPU.'')
We will show
\begin{thm}\label{CAPUprod}
{\rm (a)} The product of two APU measures is APU.
\\
{\rm (b)} The product of two CAPU measures is CAPU.
\end{thm}
Before giving the proof of Theorem \ref{CAPUprod}, we show that
it implies Theorem \ref{ConvULC}.
Recall that $\mu \in \m$ is {\em exchangeable} if $\mu(\eta)$ depends only on $|\eta|=\sum \eta_i$.
\begin{lemma}\label{CAPUexch}
For exchangeable measures, the properties ULC and CAPU are equivalent.
\end{lemma}

\mn {\em Remark.} Pemantle \cite[Theorem 2.7]{Pem} shows that
for exchangeable measures, ULC, Rayleigh, and several other
negative dependence properties coincide.
Lemma \ref{CAPUexch} adds CAPU to this list.

\mn {\em Proof of Lemma \ref{CAPUexch}.} By Corollary \ref{CAPUULC}, we need
only show that every exchangeable ULC measure is CAPU;
in fact, since conditioning
preserves both exchangeability and ULC,
it suffices to prove that an exchangeable ULC measure is APU.
But if $\mu \in \m$ is exchangeable with rank sequence
$(a_0 \dots a_n)$, then
$$\ga_i(\mu)=a_i a_{n-i} \C{n}{i}^{-1} \C{n}{n-i}^{-1}$$
so that log-concavity of (and absence of internal zeros in)
$(a_i/\C{n}{i})_{i=0}^n$
implies unimodality of $(\ga_i(\mu))_{i=0}^n$.
\qed

\mn {\em Proof of Theorem \ref{ConvULC}.} Given ULC sequences
$a=(a_0 \dots a_n)$ and $b=(b_0 \dots b_m)$,
let
$\mu \in \m_{[n]}$ and $\nu \in \m_{\{n+1 \dots  n+m\}}$ be the
corresponding
exchangeable measures; that is,
$$\mu(\eta)=\frac{a_{|\eta|}}{\C{n}{|\eta|}} ~~~ \textrm{and} ~~~ \nu(\eta)=\frac{b_{|\eta|}}{\C{m}{|\eta|}}.$$
By Lemma \ref{CAPUexch}, $\mu$ and $\nu$ are CAPU,
so that Theorem \ref{CAPUprod}(b) and Corollary \ref{CAPUULC}
give ULC for $\mu \times \nu \in \m_{[n+m]}$, completing the proof
(since the rank sequence of $\mu \times \nu$ is the convolution
of $a$ and $b$).\qed

\mn {\em Remark.} Following \cite{LiggettULC}, say an infinite nonnegative sequence $(a_0, a_1, \ldots )$ is ULC[$\infty$] if
there are no internal zeros and
$a_i^2 \geq \frac{i+1}{i} a_{i-1} a_{i+1}$ for $i \geq 1$.
The proof of Theorem \ref{ConvULC} given in \cite{LiggettULC}
allows one or both sequences to be ULC[$\infty$],
but an easy limiting argument suffices
to get this more general statement from the finite version proved here.

\mn {\em Proof of Theorem \ref{CAPUprod}.} Notice that (a) implies (b), since any measure gotten from $\mu \times \nu$ by conditioning
is the product of measures obtained from $\mu$  and $\nu$ by conditioning.

Call a nonnegative sequence $(p_0 \dots p_s)$ {\em symmetric} if $p_i=p_{s-i}$ for $i=0 \dots s$ and {\em ultra-unimodal} if
$(p_i/\C{s}{i})_{i=0}^s$ is unimodal.
Let $\mu \in \m_{[n]}$ and $\nu \in \m_{\{n+1 \dots n+m\}}$ be APU.
Then $(\C{n}{i}\ga_i(\mu))_{i=0}^n$ and $(\C{m}{i}\ga_i(\nu))_{i=0}^m$ are
symmetric and ultra-unimodal, and we want to say that
their convolution,
$(\C{n+m}{k}\ga_k(\mu \times \nu))_{k=0}^{n+m}$ is ultra-unimodal.
So we will be done if we show
\begin{lemma}\label{SymmUU}
The convolution of two symmetric ultra-unimodal sequences is ultra-unimodal
\end{lemma}
(and symmetric). It's easy to see that Lemma \ref{SymmUU}
is not true without the
symmetry assumption.

\mn {\em Proof of Lemma \ref{SymmUU}.}
Since every symmetric ultra-unimodal sequence $(p_0 \dots p_s)$
is a
positive
linear combination of sequences of the form $(\C{s}{i} {\bf 1}_{\{k \leq i \leq s-k\}})_{i=0}^s$ (and since convolution is bilinear),
it suffices to prove
that the convolution of
$(\C{s}{i} {\bf 1}_{\{k \leq i \leq s-k\}})_{i=0}^s$ and $(\C{t}{i} {\bf 1}_{\{l \leq i \leq t-l\}})_{i=0}^t$ is ultra-unimodal for all
$k,l,s,t$ with $k \leq s/2$ and $l \leq t/2$.
(Of course this is also implied by Theorem \ref{ConvULC}.)

To see this set (for $k,l,s,t$ as above)
$$f_j=\C{s+t}{j}^{-1} \sum_i \C{s}{i} {\bf 1}_{\{k \leq i \leq s-k\}} \C{t}{j-i} {\bf 1}_{\{l \leq j-i \leq t-l\}};$$
so we should show
\beq{fj}
f_j \leq f_{j+1} ~~\textrm{for all}~~ j<(s+t)/2.
\enq

\medskip
It's convenient to work with the natural interpretation
of $f_j$ as a probability.
Let $S$ and $T$ be disjoint sets with $|S|= s$ and $|T|=t$, and
let
$$Q=\{Z\sub S\cup T:
k \leq |Z \cap S| \leq s-k, ~l \leq |Z \cap T| \leq t-l\}.$$
Then $f_j=\Pr(X_j\in Q)$, where
$X_j$ is chosen uniformly from $\C{S\cup T}{j}$.
To prove (\ref{fj}),
we consider the usual coupling of $X=X_j$ and $Y=X_{j+1}$;
namely, choose
$X$ uniformly from $\C{S\cup T}{j}$
and $y$ uniformly from $(S \cup T) \sm X$,
and set $Y = X \cup \{y\}$.
We have
$$
f_{j+1}-f_j =  \Pr(X \not\in Q, Y \in Q)-\Pr(X \in Q, Y \not\in Q),$$
so should show that the right side is nonnegative.

We may assume $j\geq k+l$, since otherwise we cannot have $X\in Q$.
Then $\{X\not\in Q,Y\in Q\}$ occurs if and only if either
(i) $|X \cap S|=k-1$, $y \in S$, and $j-k+1 \leq t-l$,
or (ii) $|X \cap T|=l-1$, $y \in T$, and $ j-l+1 \leq s-k$; thus,
$$
\mbox{$\Pr(X \not\in Q, Y \in Q)
=\C{s}{k-1}\C{t}{j-k+1}\frac{s-k+1}{s+t-j}{\bf 1}_{\{j-k+1 \leq t-l\}}
 +
 \C{t}{l-1}\C{s}{j-l+1}\frac{t-l+1}{s+t-j}{\bf 1}_{\{j-l+1 \leq s-k\}}$.}
$$
Similarly (noting that $j\leq s-k+t-l$),
$\{X \in Q, Y \not\in Q\}$ occurs if and only if either
(i) $|X \cap S|=s-k$, $y \in S$, and $l \leq j-s+k$ or
(ii) $|X \cap T|=t-l$, $y \in T$, and $k \leq j-t+l $, whence
$$
\mbox{$\Pr(X \in Q, Y \not\in Q) = \C{s}{s-k}\C{t}{j-s+k}\frac{k}{s+t-j}{\bf 1}_{\{l \leq j-s+k \}}
+ \C{t}{t-l}\C{s}{j-t+l}\frac{l}{s+t-j}{\bf 1}_{\{k \leq j-t+l\}}$.}
$$
Thus, since
$$
\mbox{$\C{s}{k-1}(s-k+1)=\C{s}{s-k}k~~$ and
$~~ \C{t}{l-1}(t-l+1) = \C{t}{t-l}l$,}
$$
we will be done if we show
\beq{Ct1}
\mbox{$\C{t}{j-k+1} {\bf 1}_{\{ j-k+1 \leq t-l\}}
\geq \C{t}{j-s+k} {\bf 1}_{\{l \leq j-s+k \}}$}
\enq
and
$$
\mbox{$
\C{s}{j-l+1}{\bf 1}_{\{j-l+1 \leq s-k\}}
\geq
\C{s}{j-t+l}{\bf 1}_{\{k \leq j-t+l\}}$.}
$$

\mn
The easy verifications are similar and we just do (\ref{Ct1}):
we have
$j-k+1 \geq j-s+k$ (since $2k \leq s$) and $(j-k+1)+(j-s+k) \leq t$ (since $2j \leq s+t-1$), implying both
$
\C{t}{j-k+1}
\geq \C{t}{j-s+k}
$
and
$
 {\bf 1}_{\{ j-k+1 \leq t-l\}}
\geq  {\bf 1}_{\{l \leq j-s+k \}}
$.

\qed

\section{Consequences for Mason's conjecture} \label{section:Mason}

In this section we prove Theorem \ref{matroidCAPPsmall}.
As noted at the end of Section \ref{section:intro}, this with
Theorem \ref{CAPPULCLocal}
(or, for part (a), Theorem
\ref{CAPPULC}) immediately implies
Theorem \ref{partialMason}.
Here we do assume a (very) few matroid basics---again, 
\cite{Welshbook} and \cite{Oxley} are standard references---and
now denote matroids by $M$.
Our argument mainly follows that of
\cite{Dowling}, which, as mentioned in Section \ref{section:intro},
makes some progress on the ``LC version" of Mason's conjecture.

Given a matroid $M$ on ground set $E$, let $\Pi_{i}=\Pi_{i}(M)$ be the set of ordered partitions $(A,B)$ of $E$ with $|A|=i$ and each of $A,B$ independent.
Notice that when $|E|=2k$, APP for
uniform measure on the independent sets of $M$
is the inequality $|\Pi_{k-1}| \leq \frac{k}{k+1} |\Pi_{k}|$.

Dowling's point of departure was the observation that if
$|\Pi_k(M)| \geq |\Pi_{k-1}(M)|$
for every $k\leq t$ and every $M$ on an $E$ of size $2k$, then
for an arbitrary $M$ (on a ground set of any size) the initial portion
$(a_0 \dots a_{t+1})$ of the sequence of independence numbers is LC.
This is, of course, analogous to Theorem \ref{CAPPULCLocal}.
Note, though, that, in contrast to Theorem \ref{CAPPULCLocal}, the
implication here is quite straightforward; namely,
a natural (and standard) grouping of terms
represents the expansion of
$a_i^2\geq a_{k-1}a_{k+1}$ as a positive combination of
inequalities $|\Pi_k(M)|\geq |\Pi_{k-1}(M)|$ for various $M$'s.
(If, in analogy with (\ref{APseq}), we set
$\gb_i(\nu) = \sum\{\nu(\eta)\nu(\uone-\eta):\eta\in \Omega, |\eta|=i\}$,
then
Dowling's argument shows that $\mu\in \m$ is LC provided
each $\nu$ obtained from $\mu$ by conditioning on the values of some
$n-2k$ variables satisfies $\gb_k(\nu)\geq \gb_{k-1}(\nu)$.)

Dowling also showed
that every matroid on a ground set of size $2k\leq 14$ satisfies $|\Pi_k| \geq |\Pi_{k-1}|$ (which yields the result mentioned in
Section \ref{section:intro}).
This is mainly based on Lemma \ref{dowling2427} below and
(a version of)
the following easy observation, in which we use $d$ for degree and
``$\sim$" for adjacency.
\begin{lemma}\label{GraphWeight}
Let $G$ be a simple, bipartite graph with bipartition $X \cup Y$.
If $d{(x)} \geq 1$ for all $x \in X$ and
$\sum_{x \sim y} d(x)^{-1}
\leq C$ for all $y \in Y$,
then $|X| \leq C |Y|$.
\end{lemma}

\mn {\em Proof.}
This is standard:
$~|X|=\sum_{x \in X} \sum_{y \sim x} d(x)^{-1} =
\sum_{y \in Y} \sum_{x \sim y} d(x)^{-1} \leq C|Y|.$ \qed

\mn {\em Proof of Theorem \ref{matroidCAPPsmall}.}
Since the class of measures in question is closed under
conditioning, it's enough
to show that every matroid $M$ on a ground set $E$
of size $2k\leq 10$ satisfies
\beq{partitionAPP}
|\Pi_{k-1}(M)| \leq \frac{k}{k+1}|\Pi_k(M)|.
\enq
This is trivial when $k=1$, so
we assume $k \in \{2,3,4,5\}$.
Define bipartite graphs $G_1, G_2$
with the common bipartition $\Pi_{k-1} \cup \Pi_k$ by setting,
for $(C,D) \in \Pi_{k-1}$ and $(A,B) \in \Pi_k$,
$(C,D)\sim (A,B)$ in $G_1$ (resp. $G_2$) if $C \sub A$ (resp. $C \sub B$).
Let $G=G_1\cup G_2$.
Then, writing $r$ for rank and $d_i$ and $d$ for degrees in $G_i$ and $G$,
we have
(see \cite{Dowling}, pp. 24-27)
\begin{lemma}\label{dowling2427}
If $r(M)\geq k+2$ or $r(M)=k+1$ and $M$ has no coloops, then

\mn
{\rm (a)}
every $(A,B) \in \Pi_k$ satisfies $2 \leq d_i(A,B) \leq k$ for $i=1,2$;

\mn
{\rm (b)} every $(A,B) \in \Pi_k$ satisfies
$$\sum_{(C,D) \sim (A,B)} \frac{1}{d(C,D)} \leq \frac{1}{2}
\left( \frac{d_1(A,B)}{d_2(A,B)+1} + \frac{d_2(A,B)}{d_1(A,B)+1} \right);$$
{\rm (c)}
every $(A,B) \in \Pi_k$ with $d_1(A,B)<d_2(A,B)$ satisfies
\beq{c}
\sum_{(C,D) \sim (A,B)} \frac{1}{d(C,D)} \leq \frac{1}{2} \left( \frac{d_1(A,B)-1}{d_1(A,B)+1} + \frac{d_2(A,B)-d_1(A,B)+1}{d_1(A,B)+2} + \frac{d_1(A,B)}{d_2(A,B)+1} \right).
\enq
\end{lemma}

\mn
{\em Proof of} (\ref{partitionAPP}).
We may assume $r(M)>k$ since otherwise $\Pi_{k-1}=\0$.
Also, if $r(M)=k+1$ and $M$ has a coloop $e$,
then
$|\Pi_{k-1}(M)|=|\Pi_{k-1}(M\sm e)|$ (since every basis contains $e$)
and
$|\Pi_k(M)|=2|\Pi_{k-1}(M\sm e)|$,
so we have (\ref{partitionAPP}).

So we may assume we are in the situation of Lemma \ref{dowling2427}
(either $r(M)\geq k+2$ or $r(M)=k+1$ and $M$ has no coloops).
By Lemma \ref{GraphWeight}, it suffices to show
that for each $(A,B) \in \Pi_k$,
\beq{gweq}
\sum_{(C,D) \sim (A,B)} \frac{1}{d(C,D)} \leq \frac{k}{k+1}.
\enq
Since $d_1(A,B)=d_2(B,A)$, we may assume, using Lemma \ref{dowling2427}(a),
that $2 \leq d_1(A,B) \leq d_2(A,B) \leq k$.
If $d_1(A,B)=d_2(A,B)$, then
Lemma \ref{dowling2427}(b) bounds the left side of (\ref{gweq}) by
$d_1(A,B)/(d_1(A,B)+1)\leq k/(k+1)$.
Otherwise (i.e. if $d_1(A,B)<d_2(A,B)$),
Lemma \ref{dowling2427}(c) bounds the left side of (\ref{gweq}) by
the right side of (\ref{c}), which a little calculation shows---this is where
we use $k\leq 5$---to be at most $d_2(A,B)/(d_2(A,B)+1) \leq k/(k+1)$.\qed

\bn
{\bf Acknowledgments} We wish to thank Dave Wagner for telling us about \cite{Wagner2},
and the Isaac Newton Institute for Mathematical Sciences,
University of Cambridge, for generous support during the programme on
Combinatorics and Statistical Mechanics,
where some of this work was carried out.

\end{document}